\documentstyle[amsfonts,11pt]{article}

\newtheorem{theo}{Theorem}[section]

\newtheorem{remar}[theo]{Remark}
\newtheorem{prop}[theo]{Proposition}

\newtheorem{lemma}[theo]{Lemma}

\def\Ric{\mathop{\hbox{Ric}}}
\def\scal{\mathop{\hbox{scal}}}

\newcommand{\fdim}{\hspace*{\fill}$\Box$}
\newcommand{\dimostr}{{\bf Proof: }}

\newcommand{\R}{\Bbb{R}}

\newcommand{\complex}{\Bbb{C}}
\newcommand{\C}{\Bbb{C}}

\newcommand{\K }{K\"{a}hler}

\begin{document}

\centerline {\LARGE\bf Extremal  metrics on Hartogs domains
\footnote{The authors were supported
by the M.I.U.R. Project \lq\lq Geometric Properties of Real and
Complex Manifolds''.}}

\vspace{0.5cm}

\centerline{\small Andrea Loi and Fabio Zuddas } \centerline{\small Dipartimento di
Matematica e Informatica -- Universit\`{a} di Cagliari -- Italy}
\centerline{\small e-mail address: loi@unica.it, fzuddas@unica.it}

\vspace{0.3cm}

\begin{abstract}
An $n$-dimensional Hartogs domain $D_F$  with strongly
pseudoconvex boundary can be equipped with a natural \K\ metric
$g_F$. In this paper we prove that if $g_F$ is an extremal \K\
metric then $(D_F, g_F)$ is biholomorphically isometric to the
$n$-dimensional  complex hyperbolic space.


{\it{Subj.Class}}: 53C55, 32Q15, 32T15.

\end{abstract}

\section{Introduction and  statements of the main results}
The study of the existence and uniqueness of a  preferred \K\ metric on a given  complex manifold $M$
is a very interesting and important  area of research,
both from the mathematical and  from  the physical point of view.
Many definitions  of    canonical metrics  (Einstein, constant scalar curvature, extremal, balanced
and so on)
have been given both in the compact
and in the noncompact case
(see e.g. \cite{besse},  \cite{gau} and \cite{tian}).
In the noncompact case many important questions are still open.
For example
Yau raised the
question on the classification
of Bergman Einstein metrics
on strongly pseudoconvex domains
and S.Y. Cheng conjectured
that if the Bergman metric
on a strongly pseudoconvex domain
is Einstein, then the domain
is biholomorphic to the ball
(see \cite{si}).

\vskip 0.3cm
In this paper we are interested in extremal \K\ metrics on  noncompact manifolds.
Extremal  metrics  were introduced and christened by Calabi \cite{caex}
in the compact case as the solution for the variational problem in  a \K\ class defined by the square integral of the
scalar curvature. Therefore they are a generalization
of constant scalar curvature metrics.
Calabi himself constructs some compact manifolds with an extremal metric which cannot admit a metric with constant scalar curvature.

Only recently extremal \K\ metrics   were rediscovered  by several  mathematicians due to
their link with the  stability of complex vector bundles (see e.g.
 \cite{burns}, \cite{chtian},   \cite{fu}, \cite{levin} and  \cite{ma}).

Obviously extremal metrics  cannot be defined in the noncompact
case as  the solutions of a variational problem involving some
integral on the manifold. Nevertheless, in the compact case, one
can give an alternative definition of these metrics using  local
coordinates (see (\ref{equazestrem}) below) which  makes sense
also in the noncompact case. In this case, the existence and
uniqueness of  such    metrics are far from being  understood. For
example, only recently \cite{cheng} (see also \cite{cheng2}),
 it has been shown the existence of a nontrivial  (namely with nonconstant scalar curvature)
extremal   and complete  \K\ metric
 in a complex one-dimensional manifold.

Our main result is the following theorem   which deals with
extremal \K\ metrics on a particular class  of strongly
pseudoconvex domains, the so called  {\em Hartogs domains} (see the next
section   for their definition and main properties).

\begin{theo}\label{mainteor}
Let $(D_F, g_F)$ be an $n$-dimensional strongly pseudoconvex Hartogs domain. Assume that
$g_F$ is an extremal \K\ metric.
Then $(D_F, g_F)$ is biholomorphically isometric to the $n$-dimensional complex hyperbolic space
$({\complex}H^{n}, g_{hyp})$, where ${\complex}H^{n}$ is the  unit ball in ${\complex}^n$
and  $g_{hyp}$ denotes the hyperbolic metric.
\end{theo}


Two remarks are in order (compare with  Cheng's conjecture above).
First,  it is worth  pointing out
that, in contrast to  the Bergman metric,  $g_F$
is defined also if the domain $D_F$ is unbounded.
Secondly,  the extremality assumption  in Theorem \ref{mainteor}
is  weaker than Einstein's condition (actually it is even weaker than the constancy
of the scalar curvature).

\vskip 0.3cm

The paper is organized as follows. In the next section we recall
the definition of Hartogs domain $(D_F, g_F)$ and we analyze the  relation
between  the pseudoconvexity of $D_F$ and the \K\ condition of $g_F$.
 We also  compute its 
Ricci  and scalar curvatures. The last section is
dedicated to the proof of  Theorem \ref{mainteor}.

\section{Strongly pseudoconvex Hartogs domains}
Let $x_0 \in \R^+ \cup \{ + \infty \}$ and let $F: [0, x_0)
\rightarrow (0, + \infty)$ be a decreasing continuous function,
smooth on $(0, x_0)$. The {\em Hartogs domain} $D_F\subset {\C}^{n}$
associated to the function $F$
is defined by
$$D_F = \{ (z_0, z_1,\dots ,z_{n-1}) \in {\C}^{n} \; | \; |z_0|^2 < x_0, \ |z_1|^2 +\cdots+ |z_{n-1}|^2 < F(|z_0|^2)
\}.$$

One can prove that  the assumption of strongly pseudoconvexity of
$D_F$ is equivalent (see Proposition \ref{Kmet} below) to the fact
that the natural $(1, 1)$-form on $D_F$  given by
\begin{equation}\label{omegaf}
\omega_F = \frac{i}{2} \partial \overline{\partial}
\log \frac{1}{F(|z_0|^2) - |z_1|^2 - \cdots-|z_{n-1}|^2}
\end{equation}
is  a \K\ form on $D_F$.
The \K\  metric  $g_F$ associated to the \K\ form $\omega_F$
is the metric we will be  dealing with in the present paper.
(Observe that for $F(x)=1-x, 0\leq x< 1$, $D_F$ equals the $n$-dimensional complex hyperbolic space
${\complex} H^n$ and $g_F$ is the hyperbolic metric).
In the $2$-dimensional case this metric has been considered in \cite{eng}
and \cite{constscal} in the framework of quantization of \K\ manifolds.
In  \cite{hartogs},  the first author studied the \K\ immersions of
$(D_F, g_F)$
into finite or infinite dimensional complex space forms and \cite{alcu}  is concerned with
the existence of global  symplectic coordinates on $(D_F, \omega_F)$.
\begin{prop}\label{Kmet}
Let $D_F$ be a Hartogs domain in ${\complex}^{n}$. Then the
following conditions are equivalent:
\begin{itemize}
\item [(i)] the $(1, 1)$-form $\omega_F$  given by (\ref{omegaf})
is a \K\ form; \item [(ii)] the function $- \frac{x F'(x)}{F(x)}$
is strictly increasing, namely $-( \frac{x F'(x)}{F(x)})^{'} >0$
for every $x \in [0, x_0)$; \item [(iii)] the boundary of $D_F$ is
strongly pseudoconvex at all $z = (z_0, z_1,\dots,z_{n-1})$ with
$|z_0|^2 < x_0$ ; \item [(iv)] $D_F$ is strongly pseudoconvex.
\end{itemize}
\end{prop}
 \dimostr
 $(i)\Leftrightarrow (ii)$
 Set
 \begin{equation}\label{eqA}
 A = F(|z_0|^2)
- |z_1|^2 - \cdots - |z_{n-1}|^2.
\end{equation}
 Then $\omega_F$ is a \K\ form if and only if
the real-valued function $\Phi = -\log A$ is strictly
plurisubharmonic, i.e. the matrix $g_{\alpha\bar
\beta}=(\frac{\partial^2 \Phi}{\partial z_{\alpha}\partial \bar
z_{\beta}})$, $\alpha,\beta=0,\dots,n-~1$ is positive definite,
where
\begin{equation}\label{kform}
\omega_F=\frac{i}{2}\sum_{\alpha, \beta=0}^{n-1}g_{\alpha\bar \beta}dz_{\alpha}\wedge d\bar z_{\beta}.
\end{equation}
A straightforward computation gives
$$\frac{\partial^2 \Phi}{\partial z_0 \partial \bar z_0} =
\frac{F'^2(|z_0|^2) |z_0|^2 - ( F''(|z_0|^2) |z_0|^2 +
F'(|z_0|^2))A}{A^2} , $$
$$\frac{\partial^2 \Phi}{\partial z_0 \partial \bar z_\beta} =
- \frac{F'(|z_0|^2) \bar{z_0} z_{\beta}}{A^2}, \ \ \beta=1,\dots,n-1$$
and
$$\frac{\partial^2 \Phi}{\partial z_{\alpha} \partial \bar z_{\beta}} =
\frac{\delta_{\alpha \beta} A + \bar{z_{\alpha}} z_{\beta}}{A^2},
\ \ \alpha,\beta=1,\dots,n-1 .$$

Then, by setting
\begin{equation}\label{eqc}
C = F'^2(|z_0|^2) |z_0|^2 - ( F''(|z_0|^2) |z_0|^2 +
F'(|z_0|^2))A,
\end{equation}
one sees that the matrix $h=(g_{\alpha\bar
\beta})=(\frac{\partial^2 \Phi}{\partial z_{\alpha} \partial \bar
z_{\beta}})_{\alpha, \beta=0,\dots ,n-1}$  is given by:

\begin{small}
\begin{equation}\label{matrixg}
h= \frac{1}{A^2} \left( \begin{array}{cccccc}
C & - F' \bar z_0 z_1 & \dots  & - F' \bar z_0 z_{\alpha}&\dots & - F' \bar z_0 z_{n-1}\\
- F' z_0 \bar z_1  & A + |z_1|^2 &\dots  & \bar z_1 z_{\alpha} &\dots& \bar z_1 z_{n-1}\\
\vdots&\vdots &&\vdots&&\vdots \\
-F' z_0 \bar z_{\alpha} & z_1 \bar z_{\alpha} &\dots & A +
|z_{\alpha}|^2 &\dots&
 \bar z_{\alpha}z_{n-1}\\
\vdots&\vdots &&\vdots&&\vdots \\
-F' z_0 \bar z_{n-1} & z_1 \bar z_{n-1} &\dots & z_{\alpha} \bar
z_{n-1} &\dots & A + |z_{n-1}|^2
\end{array} \right) .
\end{equation}
\end{small}

First notice that the $(n-1)\times (n-1)$  matrix obtained by
deleting the first row and the first  column of $h$ is positive
definite. Indeed it is not hard to see that, for all $1\leq
\alpha\leq n-1$,
$$\det \left( \begin{array}{cccc} A + |z_{\alpha}|^2 & \bar
z_{\alpha} z_{\alpha+1} &\dots &
\bar z_{\alpha} z_{n-1}\\
 \vdots&\vdots &&\vdots \\
\bar z_{n-1} z_{\alpha} & \bar z_{n-1} z_{\alpha+1} &\dots& A +
|z_{n-1}|^2
\end{array} \right)=$$
\begin{equation}\label{minori principali}
 =A^{n- \alpha} + A^{n- \alpha
-1}(|z_{\alpha}|^2+\cdots+|z_{n-1}|^2) >0 .
\end{equation}

On the other hand, by the  Laplace expansion along the first row,
we get

$$\det(h) = \frac{C}{A^{2n}} [A^{n-1} + A^{n-2}(|z_1|^2+\cdots+|z_{n-1}|^2)] + $$

$$ + \frac{F' \bar z_0 z_1}{A^{2n}} \det \left( \begin{array}{cccc} -F' z_0 \bar z_1 & z_2 \bar z_1 &\dots&
z_{n-1} \bar z_{1}\\ -F' z_0 \bar z_2 & A + |z_2|^2 &\dots & z_{n-1} \bar
z_2\\
\vdots&\vdots &&\vdots \\
 -F' z_0 \bar z_{n-1} & z_2 \bar z_{n-1} &\dots & A + |z_{n-1}|^2
\end{array} \right) +\cdots + $$

$$+ (-1)^{n} \frac{F' \bar z_0 z_{n-1}}{A^{2n}} \det \left( \begin{array}{cccc} -F' z_0 \bar z_1 & A + |z_1|^2 &\dots&
z_{n-2} \bar z_{1}\\ -F' z_0 \bar z_2 & z_1 \bar z_2 &\dots& z_{n-2} \bar z_2\\
\vdots&\vdots &&\vdots \\
-F' z_0 \bar z_{n-1} & z_1 \bar z_{n-1} &\dots & z_{n-2} \bar
z_{n-1}
\end{array} \right) = $$

$$ = \frac{C}{A^{2n}} [A^{n-1} + A^{n-2}(|z_1|^2+\cdots +|z_{n-1}|^2)] + $$

$$ + \frac{F'^2 |z_0|^2 |z_1|^2}{A^{2n}} \det \left( \begin{array}{cccc} - 1 & z_2 &\dots&
z_{n-1}\\ - \bar z_2 & A + |z_2|^2 &\dots& z_{n-1} \bar z_2\\
\vdots&\vdots &&\vdots \\
- \bar z_{n-1} & z_2 \bar z_{n-1} &\dots& A + |z_{n-1}|^2
\end{array} \right) +\cdots+ $$

$$+ (-1)^{n} \frac{F'^2 |z_0|^2 |z_{n-1}|^2}{A^{2n}}
\det \left( \begin{array}{cccc} - \bar z_1 & A + |z_1|^2 &\dots&
z_{n-2} \bar z_{1}\\ - \bar z_2 & z_1 \bar z_2 &\dots& z_{n-2} \bar z_2\\
\vdots&\vdots &&\vdots \\
- 1 & z_1 &\dots& z_{n-2}
\end{array} \right) = $$

$$\frac{1}{A^{n+2}}[ C A + (C - F'^2 |z_0|^2 )(
|z_1|^2+\cdots +|z_{n-1}|^2)].$$

By substituting  (\ref{eqA}) and (\ref{eqc}) into this last
equality one gets
\begin{equation}\label{determinante metrica}
\det(h) = - \frac{F^2}{A^{n+1}} \left( \frac{x F'}{F}
\right)'|_{x=|z_0|^2} .
\end{equation}

Hence, by (\ref{minori principali}) and (\ref{determinante
metrica}),  the matrix $(\frac{\partial^2 \Phi}{\partial z_{\alpha}
\partial \bar z_{\beta}})$ is positive definite if and only if
$\left( \frac{x F'}{F} \right)'<0$.
\vskip 0.3cm

Before proving equivalence
$(ii)\Leftrightarrow (iii)$ we briefly recall some facts on complex domains
(see e.g. \cite{gurossi}).
Let  $\Omega \subseteq {\C}^{n}$
 be any complex domain of $ {\C}^{n}$   with smooth boundary $\partial\Omega$, and let $z \in \partial \Omega$.
 Assume that  there exists a smooth function
 $\rho : {\C}^n \rightarrow
{\R}$  (called  {\em defining function for $\Omega$ at $z$})
satisfying the following: for some neighbourhood $U$ of $z$, $\rho
< 0$ on $U \cap \Omega$, $\rho
>0$ on $U \setminus \overline \Omega$ and $\rho = 0$ on
$U \cap \partial \Omega$; $grad \ \rho \neq 0 $ on $\partial
\Omega$. In this case $\partial \Omega$ is said to be {\em
strongly pseudoconvex at $z$} if the {\it Levi form}

$$L(\rho, z)(X) = \sum_{\alpha,\beta=0}^{n-1} \frac{\partial^2 \rho}{\partial z_{\alpha} \partial \bar z_{\beta}}(z) X_{\alpha}\bar X_{\beta} $$
is positive definite on
$$S_{\rho}=\{ (X_0,\dots,X_{n-1}) \in {\C}^{n} \ |
\ \sum_{\alpha=0}^{n-1} \frac{\partial \rho}{\partial z_{\alpha}}(z) X_{\alpha} = 0 \}$$
(it is easily seen that this definition does not depend on the particular  defining function $\rho$).

\vskip 0.1cm

\noindent
$(ii)\Leftrightarrow (iii)$ Let now  $\Omega = D_F$ and
let us  fix   $z = (z_0, z_1,\dots,z_{n-1}) \in \partial D_F$ with  $|z_0|^2 < x_0$. 
Then,   $|z_1|^2
+\cdots +|z_{n-1}|^2 = F(|z_0|^2)$. 
In this case 
$$\rho(z_0,z_1,\dots ,z_{n-1}) = |z_1|^2 +\cdots +|z_{n-1}|^2 - F(|z_0|^2)$$
 is a (globally)
defining function for $D_F$ at $z$, the Levi form for $D_F$ reads
as
\begin{equation}\label{levi1}
L(\rho, z)(X) = |X_1|^2 + \cdots + |X_{n-1}|^2 - (F' + F'' |z_0|^2)
|X_0|^2 
\end{equation}
and
\begin{equation}\label{levi2}
S_{\rho}=\{ (X_0, X_1,\dots ,X_{n-1}) \in {\C}^{n} \ | -F' \bar z_0 X_0 +
\bar z_1 X_1 +\cdots + \bar z_{n-1} X_{n-1} = 0 \}.
\end{equation}

We distinguish two cases: $z_0=0$ and $z_0\neq 0$.
At  $z_0 = 0$ the Levi form reads as
$$L(\rho, z)(X) = |X_1|^2 +\cdots + |X_{n-1}|^2 - F'(0) |X_0|^2 $$
which is strictly positive for {\em any}  non-zero vector $(X_0,
X_1,\dots , X_{n-1})$ (not necessarily in $S_{\rho}$)
 because $F$ is assumed to be decreasing.
 
 If $z_0\neq 0$  by (\ref{levi2}) we obtain
 $X_0=  \frac{\bar z_1 X_1 +\cdots + \bar z_{n-1} X_{n-1}}
 {F' \bar z_0}$
 which,  substituted in (\ref{levi1}),  gives:
\begin{equation}\label{Levi ristretto}
L(X, z) = |X_1|^2 +\cdots + |X_{n-1}|^2 - \frac{F' + F'' |z_0|^2}{F'^2
|z_0|^2} |\bar z_1 X_1 +\cdots + \bar z_{n-1} X_{n-1} |^2.
\end{equation}

Therefore we are reduced  to show that:

\vskip 0.1cm

{\em 
$(xF' / F)' <0$ for $x \in (0, x_0)$
if and only if 
$L(X, z)$ is strictly positive
for every $(X_1,\dots ,X_{n-1}) \neq (0,\dots
, 0)$ and every $(z_0,z_1,\dots,z_{n-1}) \in \partial D_F$, $0 <
|z_0|^2 < x_0$.}

\vskip 0.3cm

If $(x F' / F)' <0$ then $(F' + x F'')F < x F'^2$ and,  since
$F(|z_0|^2) = |z_1|^2 +\cdots + |z_{n-1}|^2$, we get:
$$L(X, z) > |X_1|^2 +\cdots +
|X_{n-1}|^2 - \frac{1}{F(|z_0|^2)} |\bar z_1 X_1 +\cdots + \bar z_{n-1} X_{n-1}
|^2=$$
$$= \frac{(|X_1|^2 + \cdots + |X_{n-1}|^2)(|z_1|^2 + \cdots  + |z_{n-1}|^2) - |\bar z_1 X_1
+\cdots + \bar z_{n-1} X_{n-1} |^2}{|z_1|^2 + \cdots + |z_{n-1}|^2}$$ and the conclusion follows
by the Cauchy-Schwarz inequality.

Conversely, assume that $L(X, z)$ is  strictly positive
for every $(X_1,\dots ,X_{n-1}) \neq (0,\dots ,
0)$ and  each $z = (z_0,z_1,\dots ,z_{n-1})$
such that
$F(|z_0|^2) = |z_1|^2 +\cdots +
|z_{n-1}|^2$.
By inserting  $(X_1,\dots ,X_{n-1}) =
(z_1,\dots ,z_{n-1})$ in (\ref{Levi ristretto}) we get
$$L(z, z) = F(|z_0|^2) \left( 1 -
\frac{F' + F'' |z_0|^2}{F'^2 |z_0|^2} F(|z_0|^2) \right) >0$$
which implies $(xF' / F)' <0$.

\vskip 0.3cm

Finally, the proof of the equivalence
(ii)$\Leftrightarrow$(iv) is completely analogous to that given in
\cite{eng} (Proposition 3.4 and Proposition 3.6 ) for  the
$2$-dimensional case, to which the reader is referred. \fdim

\begin{remar}{\rm
Notice that the previous proposition  is a generalization   of
Proposition 3.6 in \cite{eng} proved there for the $2$-dimensional
case.}
\end{remar}

\vskip 0.3cm

Recall (see e.g.
\cite{kn}) that the Ricci curvature and the scalar curvature of a
\K\ metric $g$ on an $n$-dimensional complex manifold  $(M, g)$
are given
respectively by
\begin{equation}\label{riccicurvature}
{\Ric}_{\alpha \bar \beta} = - \frac{\partial^2}{\partial
z_{\alpha} \partial \bar z_{\beta}}(\log \det(h)), \ \ \alpha,
\beta = 0, \dots, n-1
\end{equation}
and
\begin{equation}\label{scalcurvature}
{\scal}_g= \sum_{\alpha, \beta = 0}^{n-1} g^{\beta \bar \alpha} {\Ric}_{\alpha \bar
\beta},
\end{equation}
where $g^{\beta \bar \alpha}$ are the entries  of the inverse
of  $(g_{\alpha\bar \beta})$, namely
$\sum_{\alpha =0}^{n-1}g^{\beta \bar \alpha}g_{\alpha\bar \gamma}=\delta_{\beta\gamma}$.

When $(M, g)=(D_F, g_F)$,
using  (\ref{matrixg})  it 
is not
hard  to  check the validity of the following equalities.

\begin{equation}\label{g00}
g^{0 \bar 0} =\frac{A}{B} F,
\end{equation}

\begin{equation}\label{gb0}
g^{\beta \bar 0} =\frac{A}{B} F' z_0 \bar z_{\beta}, \ \ \ \
\beta = 1,\dots ,n-1,
\end{equation}

\begin{equation}\label{gba}
g^{\beta \bar \alpha} =\frac{A}{B} (F' + F'' |z_0|^2) z_{\alpha}
\bar z_{\beta}, \ \ \ \ \alpha \neq \beta, \ \alpha, \beta =
1,\dots ,n-1,
\end{equation}

\begin{equation}\label{gbb}
g^{\beta \bar \beta} =\frac{A}{B} [B + (F' + F'' |z_0|^2)
|z_{\beta}|^2], \ \ \ \ \beta = 1,\dots ,n-1,
\end{equation}

where 
$$B = B(|z_0|^2)= F'^2 |z_0|^2 - F(F' + F'' |z_0|^2).$$ 
Now, set
$$L(x) = \frac{d}{dx} [x \frac{d}{dx} \log (xF'^2 -
F(F'+F''x))].$$
A straightforward computation using (\ref{determinante metrica})
and (\ref{riccicurvature}) gives:

\begin{equation}\label{ricci00}
{\Ric}_{0 \bar 0} = -L(|z_0|^2) - (n+1) g_{0 \bar 0},
\end{equation}

\begin{equation}\label{ricciAB}
{\Ric}_{\alpha \bar \beta} = - (n+1) g_{\alpha \bar \beta}, \ \ \
\alpha >0.
\end{equation}

Then, by (\ref{scalcurvature}),  the scalar curvature
of the metric $g_F$ equals
$${\scal}_{g_F} =- L(|z_0|^2) g^{0 \bar 0} -(n+1) \sum_{\alpha,\beta=0}^{n-1}
g^{\beta \bar\alpha} g_{\alpha \bar \beta} = - L(|z_0|^2) g^{0 \bar 0} -n(n+1) ,$$
which  by (\ref{g00}) reads as
\begin{equation}\label{scalgF}
{\scal}_{g_F}=  - \frac{A}{B} F L -n(n+1).
\end{equation}

 \section{Proof of the main result}\label{final}
In order to prove Theorem \ref{mainteor}, we need Lemma
\ref{mainlemma} below, interesting on its own sake,  which is a
generalization of a result  proved by the first author 
 for $2$-dimensional Hartogs domains
(see Theorem 4.8 in \cite{constscal}).

We first recall the definition of generalized scalar curvatures.
Given a  \K\ metric $g$ on   an  $n$-dimensional  complex manifold
$M$,  its {\em generalized scalar curvatures}
are the $n$ smooth functions
$\rho_0$,\dots ,$\rho_{n-1}$ on $M$
 satisfying  the following equation:
\begin{equation}\label{genscal}
\frac{\det(g_{\alpha \bar \beta} + t Ric_{\alpha \bar \beta})}{\det(g_{\alpha \bar \beta})} = 1 + \sum_{k=0}^{n-1} \rho_k t^{k+1} ,
\end{equation}
where $g_{\alpha\bar\beta}$ are the entries of the   metric in local coordinates.
Observe that for $k=0$ we recover the value of the scalar
curvature, namely
\begin{equation}\label{roscal}
\rho_0={\scal}_g.
\end{equation}
 The introduction and the study  of these curvatures  (in the compact case)
is due to K. Ogiue \cite{ogiue} to whom the reader is referred for
further results.  In particular,  in a joint paper with B.Y. Chen
\cite{CO}, he   studies  the constancy of one of the generalized
scalar curvatures. Their main result is that, under suitable
cohomological conditions,  the constancy of one of the
${\rho_k}'s, k=0, \dots ,n-1, $ implies that the metric $g$ is
Einstein.

\begin{lemma}\label{mainlemma}
Let $(D_F, g_F)$ be an $n$-dimensional  Hartogs domain. Assume
that one of its generalized scalar curvatures is constant. Then
$(D_F, g_F)$ is biholomorphically isometric to the $n$-dimensional
hyperbolic space.
\end{lemma}
\dimostr
By  (\ref{ricci00}),
(\ref{ricciAB}) we get

$$\frac{det(g_{\alpha \bar \beta} + t Ric_{\alpha \bar
\beta})}{det(g_{\alpha \bar \beta})} = (1-(n+1)t)^{n} - t L
(1-(n+1)t)^{n-1} \frac{A F}{B} .$$

 So the generalized curvatures of $(D_F, g_F)$
 are given by

\begin{equation}
\rho_{k} = (n+1)^{k}(-1)^{k+1} {n-1 \choose k}\left[
\frac{n(n+1)}{k+1} + \frac{AFL}{B} \right], \ \ \ k = 0,\dots,n-1
\end{equation}

Notice that,  for $k=0$,  we get $\rho_0= -\frac{AFL}{B}-n(n+1)=\scal_{g_F}$,
(compare with (\ref{scalgF})) in accordance with (\ref{roscal}).

Thus, $\rho_k$ is constant for some (equivalently, for any) $k
=0,\dots,n-1$ if and only if $\frac{AFL}{B}$ is constant. Since $A
= F(|z_0|^2) - |z_1|^2 - \cdots - |z_{n-1}|^2$ depends on
$z_1,\dots,z_{n-1}$ while $\frac{LF}{B}$ depends only on $z_0$,
this implies that $L =0$, i.e.

$$\frac{d}{dx} \left[x \frac{d}{dx} \log (xF'^2 - F(F'+F''x))\right]_{x = |z_0|^2} \equiv 0 .$$

Now, we continue as in the proof of Theorem 4.8 in
\cite{constscal} and conclude that $F (x)= c_1 - c_2 x, \ x=|z_0|^2$, with $c_1,
c_2 >0$, which implies that $D_F$ is biholomorphically isometric
to the hyperbolic space ${\C}H^{n}$ via the map

$$\phi: D_F \rightarrow {\C}H^{n}, \ (z_0, z_1,\dots,z_{n-1}) \mapsto \left( \frac{z_0}{\sqrt{c_1/c_2}}, \frac{z_1}{\sqrt{c_1}},\dots,\frac{z_{n-1}}{\sqrt{c_1}} \right) .$$
\fdim

\vskip 0.3cm

\noindent 
{\bf Proof of   Theorem \ref{mainteor}}
 \noindent 
 The system of  PDE's which has to be satisfied by  an extremal \K\
metric is the following (see \cite{caex}):
\begin{equation}\label{equazestrem}
\frac{\partial}{\partial \bar z_{\gamma}} \left( \sum_{\beta =
0}^{n-1} g^{\beta \bar \alpha} \frac{\partial \scal_g}{\partial
\bar z_{\beta}} \right) = 0, 
\end{equation}

for every $\alpha, \gamma = 0,\dots,n-1$ (indeed, this
is equivalent to the requirement that the (1,0)-part of the
Hamiltonian vector field associated to the scalar curvature is
holomorphic).

In order to use  equations (\ref{equazestrem}) for $(D_F, g_F)$
we write   the entries  $g^{\beta \bar \alpha}$
by separating the terms depending  only on  $z_0$ from the other terms.   More precisely, (\ref{g00}), (\ref{gb0}), (\ref{gba}) and (\ref{gbb})  can be written as follows.

$$g^{0 \bar 0} = P_{00} + Q_{00}(|z_1|^2+ \cdots +|z_{n-1}|^2), $$

$$g^{0 \bar \alpha} = \bar z_0 z_{\alpha} [P_{0a} + Q_{0a}(|z_1|^2+\cdots +|z_{n-1}|^2)], \ \ \alpha=1,\dots ,n-1, $$

$$g^{\alpha \bar \alpha} = F + P_{aa}|z_{\alpha}|^2 - (1 + Q_{aa}|z_{\alpha}|^2)\sum_{k \neq \alpha} |z_k|^2 - R_{aa}|z_{\alpha}|^4, \ \ \alpha=1,\dots ,n-1, $$

$$g^{\beta \bar \alpha} = \bar z_{\beta} z_{\alpha}[P_{ab} + Q_{ab}(|z_1|^2+\cdots +|z_{n-1}|^2)],\ \alpha \neq \beta, \ \alpha, \beta=1,\dots ,n-1, $$

where

\label{definizioniPQR}

$$P_{00} = \frac{F^2}{B}, \ \ \ Q_{00} = - \frac{F}{B}, $$

$$P_{0a} = \frac{F' F}{B}, \ \ \ Q_{0a} = - \frac{F'}{B},  $$

$$P_{aa} = \frac{F(F' + F'' |z_0|^2)}{B} -1 , \ \ \ Q_{aa} = R_{aa} =  \frac{F' + F'' |z_0|^2}{B} , $$

$$P_{ab} = \frac{F(F' + F'' |z_0|^2)}{B} , \ \ \ Q_{ab} =  - \frac{F' + F'' |z_0|^2}{B} $$

are all functions depending only on $|z_0|^2$.

We also have (cfr. (\ref{scalgF}))

\begin{equation}\label{scalG}
{\scal}_{g_F} = -n(n+1) + G (F - |z_1|^2 -\cdots -|z_{n-1}|^2)
\end{equation}

where 
$$G=G(|z_0|^2) = - \frac{L(|z_0|^2) F(|z_0|^2)}{B(|z_0|^2)}.$$

\vskip 0.3cm

Assume that $g_F$ is an extremal metric, namely equation
(\ref{equazestrem}) is satisfied. We are going to show that
$\scal_{g_F}$ is constant and hence by Lemma \ref{mainlemma}
$(D_F, g_F)$ is biholomorphically isometric to $({\complex}H^{n},
g_{hyp})$. In order to do that, fix  $i \geq 1$ and let us  consider
equation (\ref{equazestrem}) when $g=g_F$ for $\alpha = 0$,
$\gamma = i$.

We have

$$\frac{\partial \scal_{g_F}}{\partial \bar z_0} = G' z_0 (F - |z_1|^2 -\cdots -|z_{n-1}|^2) + z_0 G F'$$

$$\frac{\partial \scal_{g_F}}{\partial \bar z_i} = -G z_i. $$

So, equation (\ref{equazestrem}) gives

$$\frac{\partial}{\partial \bar z_i}\left\{ \left[P_{00} + Q_{00}\sum_{k=1}^{n-1} |z_k|^2 \right]\left[ G' z_0 (F - \sum_{k=1}^{n-1} |z_k|^2) + z_0 G F' \right] \right. -$$
$$\left.- z_0 G \left[P_{0a} + Q_{0a}\sum_{k=1}^{n-1} |z_k|^2\right]
\sum_{k=1}^{n-1} |z_k|^2 \right\} = 0,$$

namely

$$Q_{00}z_i \left[ G' z_0 (F - \sum_{k=1}^{n-1} |z_k|^2) + z_0 G F' \right] - G' z_0 z_i \left[P_{00} + Q_{00}\sum_{k=1}^{n-1} |z_k|^2 \right] -$$

$$ -z_0 G Q_{0a} z_i \sum_{k=1}^{n-1} |z_k|^2 - z_0 z_i G \left[P_{0a} + Q_{0a}\sum_{k=1}^{n-1} |z_k|^2\right] = 0 $$

Deriving again  with respect to $\bar z_i$,  we get
$$ -2
Q_{00} G' z_0 z_i^2 - 2 G Q_{0a} z_0 z_i^2 = 0.$$ Let us assume
$z_0z_i \neq 0$. This  implies $Q_{00}
G' + G Q_{0a} = 0$, i.e.  $GF' + FG' = 0$ or, equivalently,
 $G = \frac{c}{F}$ for some constant $c \in {\R}$. The
proof of Theorem \ref{mainteor} will be completed by showing that
$c =0$.  In fact, in this case $G =0$ on the open and dense subset
of $D_F$ consisting of those points such that $z_0z_i \neq 0$
and therefore, by (\ref{scalG}),
$\scal_{g_F}$ is   constant on  $D_F$.
 In order to prove that $c=0$, let us now consider
equation (\ref{equazestrem}) for $\alpha = i$, $\gamma = i$.

$$\frac{\partial}{\partial \bar z_i} \left\{ \bar z_0 z_i \left[ G' z_0
( F - \sum_{k=1}^{n-1} |z_k|^2 ) + G F' z_0 \right] \left[P_{0a} +
Q_{0a}\sum_{k=1}^{n-1} |z_k|^2  \right]  - \right.$$

$$ - G z_i \left[ F + P_{aa}|z_i|^2 -
(1 + Q_{aa}|z_i|^2)\sum_{k \neq 0,i} |z_k|^2 - R_{aa}|z_i|^4
\right] -$$
$$\left.- G z_i \sum_{k \neq 0, i} |z_k|^2 \left[ P_{ab} +
Q_{ab}\sum_{k=1}^{n-1} |z_k|^2 \right] \right\}  =0 .$$

This implies

$$ -G' |z_0|^2 z_i^2 \left[P_{0a} + Q_{0a}\sum_{k=1}^{n-1} |z_k|^2 \right] + \bar z_0 z_i^2 Q_{0a} \left[ G' z_0 ( F - \sum_{k=1}^{n-1} |z_k|^2 ) + G F' z_0 \right] -  $$

$$- P_{aa} G z_i^2 + G z_i^2 Q_{aa}  \sum_{k \neq 0, i} |z_k|^2 + 2 G z_i^3 \bar z_i R_{aa} - G z_i^2 Q_{ab}  \sum_{k \neq 0, i} |z_k|^2. $$

If we divide by $z_i^2$ (we are assuming $z_i \neq 0$) and derive
again the above expression with respect to $\bar z_i$ we get
$$ -G' |z_0|^2 Q_{0a} + G R_{aa} = 0 .$$
By the definitions made at page \pageref{definizioniPQR}
this is equivalent to

$$ \frac{G' F' |z_0|^2 + G (F' + F'' |z_0|^2)}{B} = 0, $$
i.e. $(G F' x)' = 0, x=|z_0|^2$. Substituting $G = \frac{c}{F}$
in this equality we get $c (\frac{F' x}{F})' = 0$.
Since $(\frac{F' x}{F})' <0$
(by (ii) in Proposition \ref{Kmet}) $c$ is forced to be zero,
and this concludes the proof.

\small{}

\end{document}